\documentclass[9pt,twocolumn,twoside]{pnas-new}
\templatetype{pnasresearch article} % Choose template 
\title{Comparative analysis of Forman-Ricci curvature versions applied to the persistent homology of networks}

\author[1]{Sergio Serrano de Haro Iváñez}

\affil[1]{Contact: sergioserrano1998@gmail.com}

\leadauthor{Serrano de Haro Iváñez} 

\significancestatement{In the world data science it is crucial to be able to extract simple and computable measurements from unapproachably huge data sets. Persistent homology is a topological tool that achieves this by essentially assigning an appearance time to the elements of a space and studying the evolution of the homology groups as the objects appear over time. In the case of unweighted networks, Forman-Ricci curvature has been proposed as a relevant and easily computable way of assigning such appearance times. However, Forman-Ricci curvature is an infinite series, and so one must decide which terms to include in its calculation. In this article we provide an overview of both Forman-Ricci curvature and persistent homology, and we compare the persistence results obtained with some of the most common ways of computing Forman-Ricci curvature.}

\authorcontributions{Sergio Serrano de Haro Iváñez is responsible for the research and elaboration of the manuscript.}
\authordeclaration{No competing interests to declare.}

\begin{abstract}
We provide an overview of Forman-Ricci curvature and persistent homology, and how their combination can be applied to the study of networks. We discuss how the usually employed augmented Forman-Ricci curvature formula, only valid for quasiconvex augmented networks, can be extended to the non-quasiconvex case.  We apply three versions of quasiconvex Forman-Ricci curvature (plain, triangle-augmented, and pentagon-augmented) to build time filtrations on non-quasiconvex networks, both model and real-world. Our results suggest that triangle-augmented curvature should be used until the non-quasiconvex formula is further studied, as plain curvature omits too much information, and quasiconvex pentagon-augmented curvature is too rough of an approximation and significantly distorts the results.
\end{abstract}

\begin{document}

\maketitle
\thispagestyle{firststyle}
\ifthenelse{\boolean{shortarticle}}{\ifthenelse{\boolean{singlecolumn}}{\abscontentformatted}{\abscontent}}{}

% If your first paragraph (i.e. with the \dropcap) contains a list environment (quote, quotation, theorem, definition, enumerate, itemize...), the line after the list may have some extra indentation. If this is the case, add \parshape=0 to the end of the list environment.
\dropcap{W}ith network science becoming increasingly important, more and more sophisticated techniques are being developed for its study. Among them, we find two groups of approaches: the topological and the geometric ones. The first include those derived from the field of Topological Data Analysis, which has a very powerful tool in the form of persistent homology \cite{ROADMAP}. Persistent homology is a noise-resistant technique that allows one to compactly represent how data features evolve across different scales, and which has already shown its applications in the study of, for example, contagion \cite{CONTAGION}, financial \cite{FINANCIAL}, or tumor vascularity \cite{CANCERPH} networks.

On the other hand, geometric approaches include the discretisation of classical differential geometry notions. One such concept is Ricci curvature, which in its original setting measures how both volumes grow and geodesics diverge at each point of a Riemannian manifold. There are multiple adaptations of Ricci curvature to the discrete world, two of which stand out: those due to Ollivier \cite{OLLIVIER1,OLLIVIER2}, and Forman \cite{FORMAN}, each respectively capturing one of the two aspects of the smooth version. In this article we will focus on the latter, which has found multiple applications in network science, like dynamic change detection \cite{CHANGEDETECTION}, anomaly detection \cite{FORMANDETECTION}, or the study of the persistent homology of networks \cite{FORMANRICCIANDPERSISTENCE}. This last article, in which Forman-Ricci (FR from now on) curvature was proposed as a way to impose time filtrations on networks, is the starting point of our studies. FR curvature can be written as an infinite series, and so requires the choice of a truncation point to compute it. The arguably most common way of doing so was used in \cite{FORMANRICCIANDPERSISTENCE}, and here we compare the results of that truncation with other, less common (but potentially necessary) ones.

\section*{Forman-Ricci curvature}
 
As the main focus of the article is not the theoretical construction of FR curvature, we will just provide an overview of it, referring the reader to Forman's original article \cite{FORMAN} for a more detailed and formal discussion. FR curvature was originally designed for weighted CW-complexes, which is a big class of topological objects that contains undirected weighted networks (for an introduction to CW-complexes and other algebraic topology concepts, we refer the reader to the classical book by Hatcher \cite{HATCHER}). The original, general formula for FR curvature (see Theorem 2.2 in \cite{FORMAN}) requires defining and keeping track of an orientation on the whole CW-complex, which could slow computations down. Fortunately, networks without multiedges are a \textit{special} type of weighted CW-complex named ``quasiconvex'' (see Definition 0.1 in \cite{FORMAN}), which essentially means that the intersection of two $n$-dimensional objects of the complex is at most one $(n-1)$-dimensional object (for networks, any two edges must share at most one vertex -- that is, no multiedges). For quasiconvex complexes, the FR curvature formula can be simplified to one that no longer needs an orientation (see Theorem 3.10 in \cite{FORMAN}). Using this, one gets that the Forman-Ricci curvature of a network edge $e=(v_1,v_2)$ is \cite{FORMANCURVATUREFORCOMPLEXNETWORKS}

\vspace{-2.5ex}
\begin{align}\label{eq:1}
    \mathcal{F}(e)=&\ \omega(v_1)+\omega(v_2)\\&-\omega(e)\cdot \left(\sum_{e_{v_1}} \frac{\omega(v_1)}{\sqrt{\omega(e)\omega(e_{v_1})}} +\sum_{e_{v_2}} \frac{\omega(v_2)}{\sqrt{\omega(e)\omega(e_{v_2})}}\nonumber \right),
\end{align}

where $e_{v_i}$ denotes an edge that shares node $v_i$ with $e$, and $\omega(\cdot)$ denotes the weight of an element. We remark that equation [\ref{eq:1}] is valid for the broadest sense of weighted networks, where both nodes and edges can be assigned independent weights. In the case of unweighted networks, however, it is direct to check that the formula simplifies to:

\vspace{-2.5ex}
\begin{align}
    \mathcal{F}(e)&= 4 - \deg(v_1) - \deg(v_2). \label{eq:2}
\end{align}

%For networks with just weighted edges, both theoretical and experimental results indicate that assigning the weight $1$ to every node will provide a curvature measure more relevant than assigning them their degree weight \cite{SYSTEMATIC}.

From now on we will refer to [\ref{eq:2}] as ``plain'' FR curvature, as the generality of Forman's construction actually allows for a more general definition \cite{COMPARATIVE}. This generalisation, usually called ``augmented'' curvature, adds terms that account for the cycles of the network. A  network is originally a one-dimensional object (it has 0-dimensional nodes and 1-dimensional edges), and we may turn it into a two-dimensional one by inserting flat faces that fill the spaces encircled by its cycles. This new two-dimensional object is too a CW-complex, and Forman's work allows us to compute the curvature of its edges. Now, if we assume that this augmented network is still quasiconvex (that is, any two inserted faces share at most one edge), applying Forman's quasiconvex formula gives the following augmented FR curvature for an edge $e=(v_1,v_2)$:

\vspace{-2.5ex}
\begin{align}\label{eq:3}
    \mathcal{F}_\text{Aug}&
    (e)= \omega(v_1)+\omega(v_2)+\sum_{f_e} \frac{\omega(e)^2}{\omega(f_e)}
    \\ &-\omega(e)\cdot \sum_{e'\not= e} \left|\sum_{f_e > e'} \frac{\sqrt{\omega(e)\omega(e')}}{\omega(f_e)} - \sum_{v_i<e'} \frac{\omega(v_i)}{\sqrt{\omega(e)\omega(e')}}\right|\nonumber,
\end{align}

where $\sum_{e'\not= e}$ sums over all edges different from $e$, $\sum_{v_i<e'}$ over the nodes shared by $e$ and $e'$, and $f_e$ denotes the faces containing $e$, so  $\sum_{f_e > e'}$ sums over all faces containing both $e$ and $e'$.  Once again, the formula greatly simplifies for unweighted networks (see Lemma 1 of \cite{COARSE}):

\vspace{-2.5ex}
\begin{align}\label{eq:4}
    \mathcal{F}_\text{Aug}(e)&= \mathcal{F}(e) +3\cdot |f^3_e| + 2\cdot |f^4_e| + \dots + (6-n)\cdot |f^n_e| + \cdots 
\end{align}where $|f^n_e|$ denotes the number of faces with $n$ sides (i.e. faces inserted in a $n$-cycle) that contain $e$. We remark that equation [\ref{eq:4}] is only valid when no two faces (i.e. cycles) of the network share more than one edge, which does not tend to be the case in general networks. This has generated a couple of different approaches to avoid the use of orientations and apply [\ref{eq:4}]: the first is working in a context where the augmented networks \textit{will be} quasiconvex by construction \cite{QUANTUMPHYSICSFORMAN}. The second consists in working with general networks, but only inserting the faces corresponding to 3-cycles (i.e. triangles) \cite{COMPARATIVE,FORMANRICCIANDPERSISTENCE}. This ensures that the augmented network will be quasiconvex, as if two triangles share more than one edge, then they must be the same one. From now on, we will denote the curvature obtained this way by $\mathcal{F}_\triangle$. 

This last approach is also computationally practical, as counting large cycles becomes prohibitively expensive quickly, so one has to sacrifice the information they provide to achieve assumable computation times. Even if [\ref{eq:4}] truncated at the triangles has been pretty well studied so far, to our knowledge little work has been done on studying the higher order terms, both for the quasiconvex case (as higher terms may be needed to correctly approximate $\mathcal{F}_\text{Aug}$), and the non-quasiconvex one (where [\ref{eq:4}] could be a good enough approximation of the general formula). We believe this to be an interesting and potentially fruitful avenue of investigation, although its full treatment is out of the scope of this paper.

\begin{figure}
\centering
\includegraphics[width=.7\linewidth]{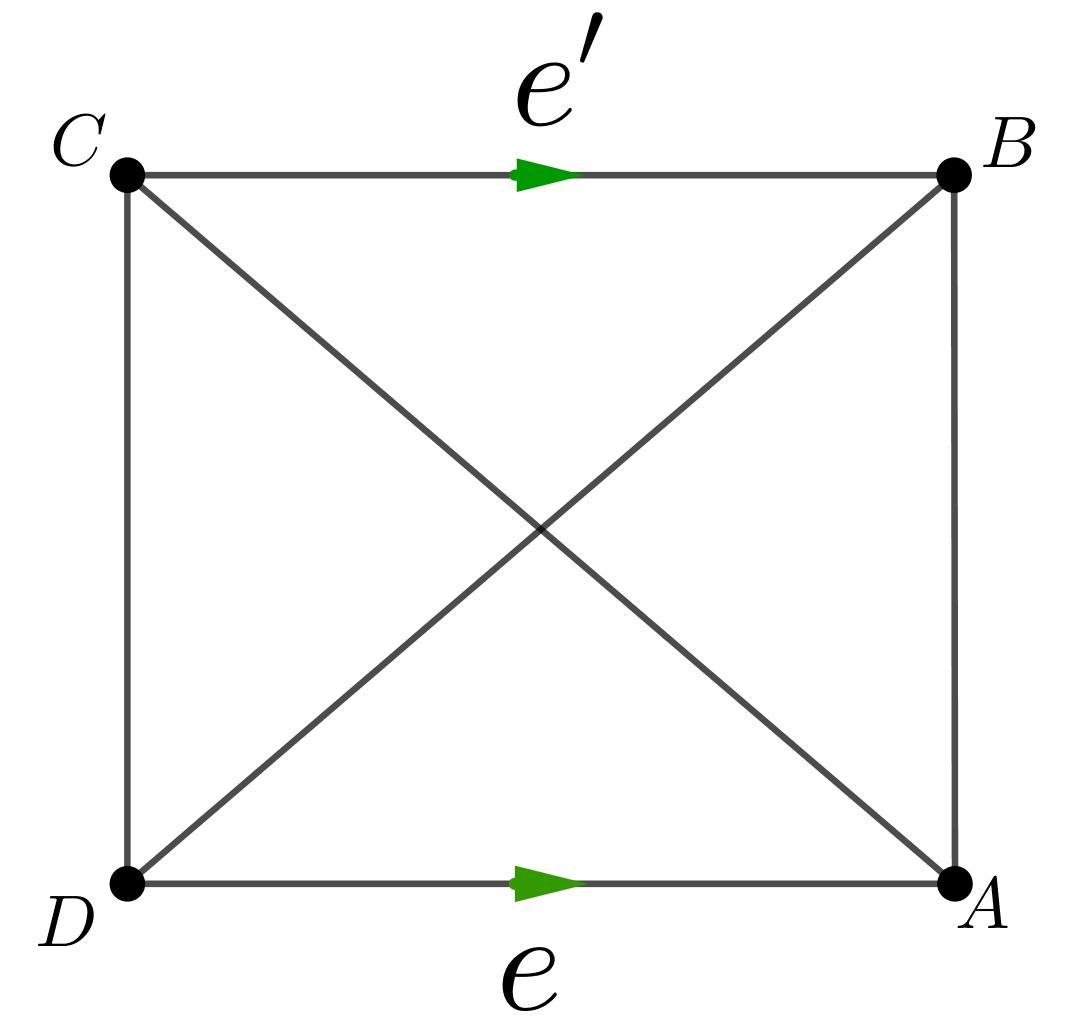}
\caption{Complete graph of 4 nodes. The green arrow marks the choice of orientation for $e$ and $e'$.}
\label{fig:K4}
\end{figure}

While we will not delve into the topic here, we would like to give a rough idea of the difference between the general and quasiconvex augmented formulas, as we will need to keep it in mind for the experimental part of the paper. Equation [\ref{eq:4}] is actually a rewrite of the following \cite{COARSE}:

\vspace{-2.5ex}
\begin{align}\label{eq:5}
    \mathcal{F}_\text{Aug}(e)&= 2 + |f_e| -\#\{e'\not= e : e'\text{ shares a node \textit{or} face with } e\},
\end{align}

where the \textit{or} is exclusive, and $|f_e|$ denotes the number of faces that contain $e$. Thus, augmented FR curvature just count faces and the edges that satisfy some condition. Now, to extend the edge count to the non-quasiconvex case we must face a multiplicity issue. For example, see the network of Figure \ref{fig:K4}: edges $e$ and $e'$ share two 4-faces (squares $ABCD$ and $ACBD$). To compute $\mathcal{F}_\text{Aug}(e)$, should we count $e'$ once (as it is only one edge), or twice (as it shares two faces with $e$)? The answer needs an orientation: Forman's general formula applied to augmented networks gives  

\vspace{-2.5ex}
\begin{align}\label{eq:6}
    \mathcal{F}&_\text{Aug}(e)= 2 + |f_e| - \sum_{e':\exists v <e,e'}\Big|\#\{f > e, e'\}-1\Big|
    \\ &-\sum_{e':\nexists v <e,e'} \Big|\#\{f > e, e' : e \xrightarrow{+}_fe'\}-\#\{f > e, e' : e  \xrightarrow{-}_f e'\}\Big|,\nonumber
\end{align}

where we respectively sum over edges sharing and not sharing a node with $e$. Both sums account for the multiplicity of overlapping faces, and the orientation comes into play in the second. Orienting $e$ gives a direction to travel through the cycles that $e$ is part of. The terms $e \xrightarrow{+}_f e'$ and $e \xrightarrow{-}_f e'$ then denote that going through the cycle that is the border of $f$ we go through $e'$
in, respectively, the same or opposite direction as its orientation. For instance, in Figure \ref{fig:K4} the orientation of $e'$ coincides with the direction of travel when going through square $ACBD$, but is opposite to the direction of travel of $ABCD$. It is direct to check that [\ref{eq:6}] reduces to [\ref{eq:5}] for quasiconvex augmented networks. In short, the general formula contemplates a \textit{signed} multiplicity through orientation which is lost when working in the quasiconvex case. As an example of the type of error one can expect by approximating [\ref{eq:6}] by [\ref{eq:4}], computing both values for $e$ in Figure \ref{fig:K4} gives respectively 2 and 8, the offset precisely caused by [\ref{eq:4}] not taking into account the overlapping faces of the network and their orientations.

Having briefly described the role of orientation in FR curvature, we turn our attention back to [\ref{eq:4}]. We would like to remark the work done in \cite{QUANTUMPHYSICSFORMAN}, where it was proven that, for some specific types of networks, using the augmented curvature just up to pentagonal faces is equivalent to measuring Ollivier-Ricci curvature. Ignoring the deeper ties between curvatures that this entails, we just want to emphasise it as evidence that the computationally necessary truncation of formula [\ref{eq:4}] does not completely spoil its relevance as a curvature measure. 

\begin{figure}
\centering
\includegraphics[width=.8\linewidth]{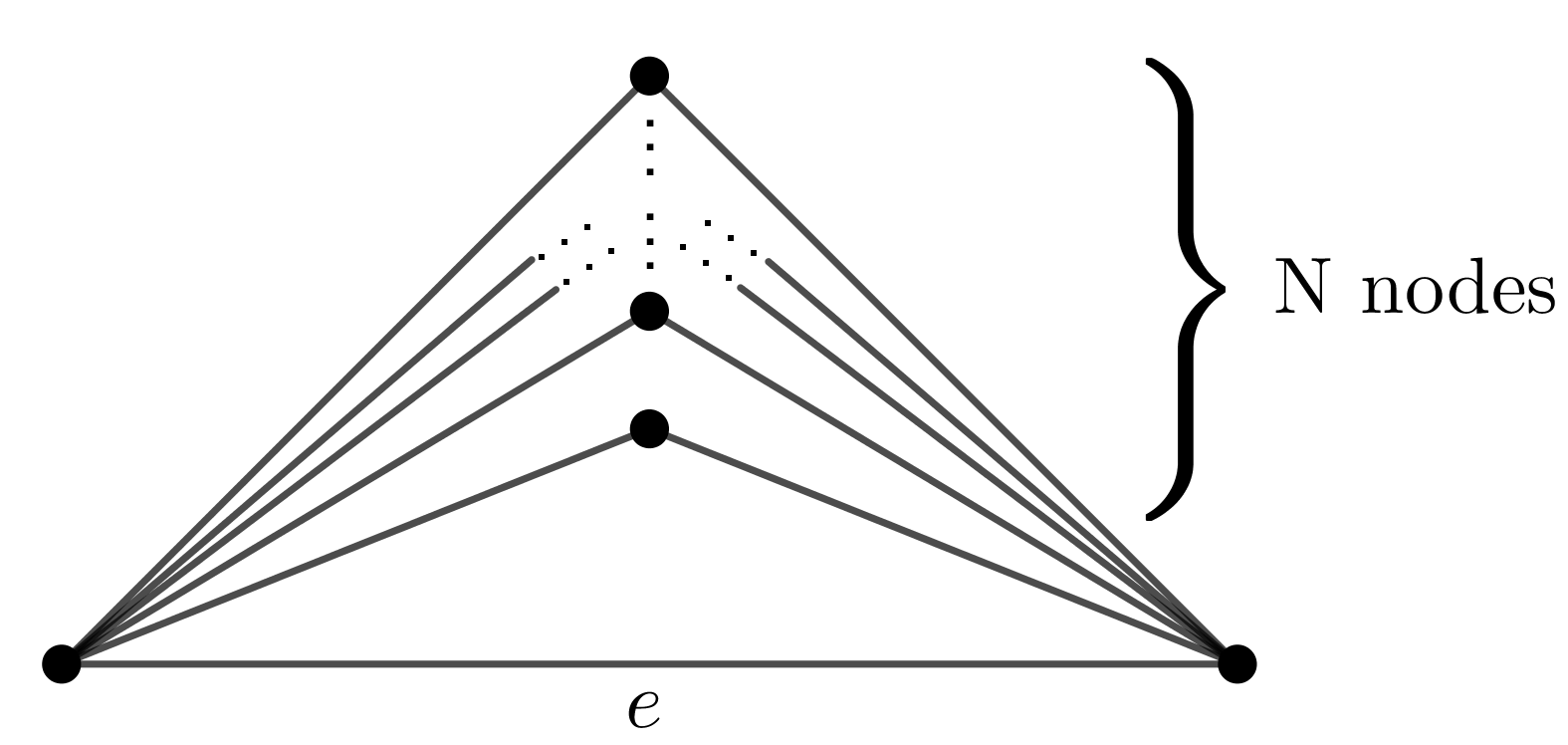}
\caption{Example of a network where an edge has very disparate plain and triangle-augmented FR curvatures.}
\label{fig:example}
\end{figure}

Now, being FR curvature a quantity with a geometric basis, it does have the potential to provide insight about network properties. This is indeed the case, as experimental studies show that both in its plain \cite{FORMANCURVATUREFORCOMPLEXNETWORKS,SYSTEMATIC}, and triangle-augmented \cite{COMPARATIVE} versions, highly negative FR curvature is generally a good indicator that an edge is important for the network's structure. However, we would like to note that even if \cite{COMPARATIVE} showed that both constructions of FR curvature seem to share this characteristic for real-world networks, it is quite easy to construct examples where their values greatly differ -- see for example the network of Figure \ref{fig:example}, where while every edge $e'\not=e$ shows little variance between the two values, $\mathcal{F}(e') = 1 - N$, $\mathcal{F}_\triangle(e') = 4 - N$, the curvatures of edge $e$ are strikingly different: $\mathcal{F}(e) = 2 - 2N$, $\mathcal{F}_\triangle(e) = 2 + N$. In a similar fashion, one can build networks such that adding any specific term of [\ref{eq:4}] greatly modifies the value of the curvature of a single edge. As we will see, this type of discordance between plain and augmented curvature is not relegated to artificial networks, but can also be found in organic ones.

%We remark that both definitions are valid only for weighted undirected networks, but can be generalised to directed networks \cite{DIRECTED} and hypernetworks \cite{HYPERNETWORKS}.

\section*{Persistent Homology}

We now present the second theoretical half of the paper. Again, we will provide a brief and intuitive overview on persistent homology, and we refer the reader to \cite{ROADMAP} for a more in-depth review of the topic.

Homology is one of the central objects of study in algebraic topology. Roughly put, the homology groups of a topological space $X$ are a sequence of vector spaces $\{H_k(X)\}_{k\geq0}$, such that $H_k(X)$ is the free vector space generated by the number of $k$-dimensional holes of $X$ -- that is, $\dim H_k(X)= \# k\text{-holes of }X$. Intuitively, 0-holes are connected components, 1-holes are ``usual'' holes, and 2-holes are cavities or voids -- a way to visualise $k$-holes is through spheres: the $n$-sphere has exactly one $n$-dimensional hole\footnote{Of course, there is much more to homology groups: for example, they can be just groups instead of vector spaces, and so they can potentially have torsion, which complicates a bit the interpretation of the elements of $H_k(X)$ as holes. However, the explanation here should be enough to keep an intuitive idea of homology.}. For networks, the 0-holes are their connected components, and the 1-holes are generated by their cycles -- however, we would like to remark that even if there is a bijection between components and 0-holes, generally there is \textit{not} one between cycles and 1-holes: holes are in truth equivalence classes, and so multiple cycles of the network can correspond to the same topological 1-hole.

Persistent homology is a tool that, given a space $X$ and an increasing sequence of nested subspaces (also known as a filtration) $X_0\subset X_1 \subset \dots \subset X_n = X$, studies how the homology groups of $X_t$ evolve as $t$ increases. We may intuitively think that the space $X_0$ discretely grows over time until it becomes $X$, and we study how different holes appear and disappear as $X_0$ develops and approaches $X$.

The usual way of representing the results of this study is through what it is known as a barcode. The $k$-barcode is a multiset of intervals $B_k=\{[b_i,d_i]\}_{i\in I}$ that registers the birth and death time of all the $k$-dimensional holes that arise in the filtration -- each interval corresponding to a different hole. It is conventional to represent holes that never disappear (that is, holes present in the final step of the filtration, the original space $X$) by intervals with death time $d_i=\infty$. Figure \ref{fig:barcodesCube} shows a standard way of displaying barcodes: each horizontal bar represents a different hole, and the $x$ coordinates of its endpoints are, respectively, its appearance and disappearance times.  One of the key points of persistent homology is that it is robust to perturbations -- that is, small modifications of either the space $X$ or the filtration will only result in equally small modifications of the barcode.

The bottleneck distance is a standard way of qualitatively comparing how close two different barcodes are. Let $B,C$ be two barcodes, $B_0\subset B, C_0\subset C$ subsets, and $\phi: B_0 \rightarrow C_0$ a bijection. Let $\varepsilon_\phi$ be the smallest $\varepsilon$ such that:

\begin{enumerate}
  \item For $[b,d]\in B_0$, if $[b',d']:= \phi\left([b,d]\right)$, then $|b'-b|\leq \varepsilon$ and $|d'-d|\leq \varepsilon$.
  \item For $[b,d]\in B_0^c \cup C_0^c$ we have $d-b \leq 2\varepsilon$.
\end{enumerate}
The bottleneck distance is then the minimum such $\varepsilon_\phi$,

\vspace{-2.5ex}
\begin{align}\label{eq:7}
    d_{\text{Bot}}(B,C)&= \min_{B_0,C_0,\phi} \varepsilon_\phi.
\end{align}

In a nutshell, one finds the best fitting pairing between the bars of $B$ and $C$, and the bottleneck distance is the greatest difference between the two bars of one of such pair.

Now, networks fall into the category of simplicial complexes, a kind of purely combinatorial objects whose homology can thus be computed combinatorically -- i.e. with a computer. Hence, networks are all set for a homological study. However, we still need a way to assign an appearance time to the elements of a network if we want to use the tools of persistent homology. For weighted networks one can use the weights for this purpose \cite{STRATA}, but unweighted networks require the extra step of deciding an appropriate time filtration. The triangle-augmented Forman-Ricci curvature has been proposed as a way of doing so \cite{FORMANRICCIANDPERSISTENCE}, with the idea of having relevant edges (at least curvature-wise, i.e. with very negative curvature) appearing soon in the filtration, thus having key structural features of the network manifesting early in the process.

However, in the original study in \cite{FORMANRICCIANDPERSISTENCE} only the triangle-augmented FR curvature was employed for the construction of time filtrations. As we have discussed, even if both plain and augmented FR curvature share some properties, one should not expect them to take the same values -- see Figure \ref{fig:example}. However, if both curvatures behave similarly enough, one could save considerable computation time by working with the plain version. On the other hand, a priori we should not assume $\mathcal{F}_\triangle$ to be an accurate enough approximation of augmented FR curvature -- maybe additional terms are needed.

Thus, we wonder how the persistence results will vary if we build time filtrations with different versions of FR curvature. Here we will study three: plain, triangle-augmented, and pentagon-augmented (that is, truncating formula [\ref{eq:4}] after the 5-cycle term), which we will denote by $\mathcal{F}_{\pentagon}$. We will compare the plain and triangle-augmented results to see the feasibility of working with plain curvature, and we will compare the two augmentations to see to what extent the higher terms influence the results. There are two reasons for the pentagon truncation choice: firstly, counting 6-cycles starts to be too time-consuming for some of the examined networks. Secondly, as mentioned, it has been shown that in some cases $\mathcal{F}_{\pentagon}$ is essentially Ollivier-Ricci curvature, so we can hope for $\mathcal{F}_{\pentagon}$ to already have a non-trivial, meaningful geometric value in the general case.

As a final note, the augmented networks we will work with most certainly will not be quasiconvex, so $\mathcal{F}_{\pentagon}$ is going to show an additional offset from the actual curvature value, which we need to account for in the analysis of our results.

\section*{Materials and methods}

We compute three different truncations of [\ref{eq:4}] for both model and real-world unweighted networks, and use it to build time filtrations on them and calculate their persistent homology. All simulations have been run in \texttt{Python}, and the packages \texttt{NetworkX} \cite{NETWORKX} and \texttt{GUDHI} \cite{GUDHI} have been used respectively for the treatment of networks and persistent homology. We have worked with the following set of models, which resembles the one reviewed in \cite{FORMANRICCIANDPERSISTENCE}:
\begin{itemize}
  \item \textbf{Erdös-Rényi (usual and bipartite)}: parameters $n\in\mathbb{N}$ and $p\in [0,1]$. Graph of $n$ nodes where each edge has an independent random probability $p$ of being present.
  \vspace{-1.3ex}
  \item \textbf{Watts-Strogatz}: parameters $n,k\in\mathbb{N}$ and $p\in [0,1]$. Ring graph of $n$ nodes where each node starts connected to its $k$ closest neighbours. Then, each edge is rewired with independent probability $p$.
  \vspace{-1.3ex}
  \item \textbf{Barabási-Albert}: parameters $n,k,n_0\in\mathbb{N}$. Graph of $n$ nodes created from an initial graph of $n_0$ nodes that is grown by successively attaching new nodes of degree $k$. The new edges are added at random, with probabilities proportional to the degrees of the existing nodes.
  \vspace{-1.3ex}
  \item \textbf{Random cube graph}: parameters $n\in\mathbb{N}$, $r\in\mathbb{R}^+$. A total of $n$ points are chosen uniformly at random in the 3-dimensional unit cube. An edge is drawn between two points if their distance is less than $r$.
\end{itemize}

All models have been set to 1000 nodes and expected average degrees of 2, 4, 6, 8, and 10. For the Watts-Strogatz model, the chance of rewiring an edge has been set to $0.5$. For the Barabási-Albert model, the initial graph is a star graph with $n_0=k$. For each graph, 30 simulations have been performed, and the bottleneck distances presented in this articles are the average of said simulations. The amount of simulations has been chosen so that one can apply the central limit theorem to find confidence intervals for the means.

Moreover, the following real-world networks have been studied:

\begin{itemize}
  \item \textbf{Western US Power Grid}: Network of the US Western States Power Grid, where nodes are transforms or power relay points, and edges are power lines between them. It has $4941$ nodes, $6594$ edges, and average degree of $\sim\hspace{-0.2em}1.3$ \cite{USPOWER}.
  \vspace{-1.3ex}
  \item \textbf{Sister cities}: Sister cities network with $14274$ nodes, $20573$ edges, and average degree of $\sim\hspace{-0.2em}1.4$ \cite{SISTERCITIES}.
  \vspace{-1.3ex}
  \item \textbf{Network science collaboration}: Co-authorship network in the area of network science, with $1461$ nodes, $2742$ edges, and average degree of $\sim\hspace{-0.2em}1.9$ \cite{NETWORKCOLLABORATION}.
  \vspace{-1.3ex}
  \item \textbf{Pretty Good Privacy}: Giant connected component of the interaction network of Pretty Good Privacy users. With $10680$ nodes, $24316$ edges, and average degree of $\sim\hspace{-0.2em}2.3$ \cite{PRETTYGOODPRIVACY}.
  \vspace{-1.3ex}
  \item \textbf{ArXiv GR-QC}: General Relativity and Quantum Cosmology arXiv collaboration network, with $5242$ nodes, $14496$ edges, and average degree of $\sim\hspace{-0.2em}2.8$ \cite{ARXIV}.
  \vspace{-1.3ex}
  \item \textbf{Email exchange}: Network of email exchanges of Universitat Rovira i Virgili members, with $1133$ nodes, $5451$ edges, and average degree of $\sim\hspace{-0.2em}4.8$ \cite{EMAIL}. 
  \vspace{-1.3ex}
  \item \textbf{Collins yeast interactome (2007)}: Network of protein-protein interactions, with $1622$ nodes, $9070$ edges, and average degree of $\sim \hspace{-0.2em}5.6$ \cite{COLLINSYEAST}.
\end{itemize}

We have tried our best to find real networks with the same average degree range as the models, but we have not been able to find any network with an average degree higher than 6 that we could compute in sensible time.

Now, the procedure for building a time filtration for each network and curvature replicates the original one in \cite{FORMANRICCIANDPERSISTENCE}, and goes as follows: for any version of FR curvature, one computes the value of all edges, and rescales them so that they lay on the interval $[0,1]$. This is done through the formula

\begin{equation}\label{eq:norm}
    \mathcal{F}_\text{N}(e) = \frac{\mathcal{F}(e)-\mathcal{F}_\text{min}}{\mathcal{F}_\text{max}-\mathcal{F}_\text{min}},
\end{equation}

where $\mathcal{F}_\text{max}$ and $\mathcal{F}_\text{min}$ denote the maximum and minimum curvatures among all the edges. The appearance times of the elements of the network are then set to be the following:
\begin{itemize}
  \item \textbf{Nodes}: minimum appearance time among incident edges. If the node has no edges, the time is set to 1.
  \vspace{-1.3ex}
  \item \textbf{Edges}: normalised curvature [\ref{eq:norm}].
  \vspace{-1.3ex}
  \item \textbf{Triangles}: maximum appearance time among forming edges.
\end{itemize}

That is, the appearance time of an edge is determined by their curvature. Nodes are forced to arise by the edges: once an edge appears, its endpoints appear too if they were not already there. Finally, a solid 2-dimensional triangle appears together with the last of its edges. The appearance time of triangles coincides with how they appear in the Vietoris-Rips complex, which is one of the most popular and easily computable ways of creating time filtrations in topological data analysis \cite{ROADMAP}. Note too that whenever any 3-cycle is created, it is instantly filled with a face, so the 1-holes will be generated by cycles of length $4$ or more. According to standard practice, the homology groups are calculated over the field $\mathbb{Z}/2$, which greatly reduces computation complexity and time.

Note that by normalising the curvatures all finite bars will be contained in the interval $[0,1]$. Thus, the greatest bottleneck distance between two barcodes is 1 (if they have the same number of infinite bars, otherwise the distance is $\infty$). Moreover, any bottleneck distance greater than $0.5$ is a result of the difference between infinite bars.

A comment on the filtration construction: note that we are inserting the solid triangles in the $\mathcal{F}$-induced filtration, although plain FR does not take them into account. This is not a problem, as the two roles the triangles play (affecting the appearance time of the edges through FR curvature and acting as a face in the computation of homology) are completely different. Moreover, if we are interested in comparing the persistence results of different versions of FR curvature, we should always calculate the homology of the same object --- removing the triangles in plain FR curvature would give us a very different complex with a potentially much bigger 1st homology group (i.e. with many more 1-holes), and an always null 2nd homology group (as there would not be any faces to create 2-holes). For the same reasons, not inserting the solid squares and pentagons in the filtration while working with $\mathcal{F}_{\pentagon}$ is not an issue, but necessary for our study. We could potentially have inserted them in the filtrations of all curvatures, but the software we have used (which works with simplicial complexes) can only insert triangular faces. This should however not be a concern, as we are interested in the difference between time filtrations on the same object, not in the object itself. As a final remark, note that inserting squares and pentagons would not affect the behaviour of the 0th homology group (i.e. the connected components), so the study of the 0-barcodes would be the exact same.

\section*{Results and discussion}

\begin{figure}%[tbhp]
	\centering
	\begin{subfigure}[b]{0.47\linewidth}
		\centering
		\includegraphics[width=1\linewidth]{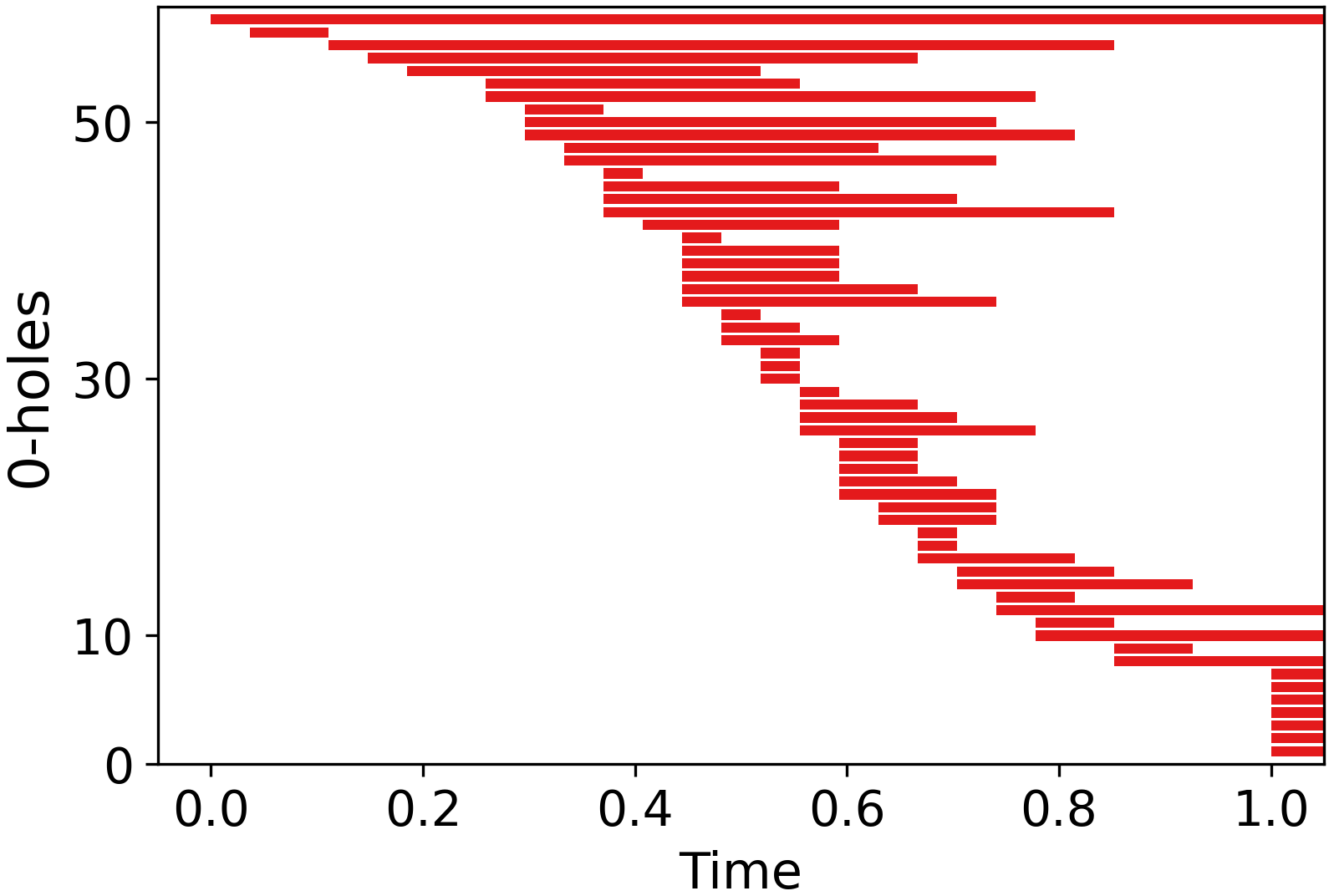} 
		\caption{$\mathcal{F}$-induced filtration.} 
		\label{fig:barcodesCube_a}
		\vspace{4ex}
	\end{subfigure}%
	\hspace{1em}%
	\begin{subfigure}[b]{0.4775\linewidth}
		\centering
		\includegraphics[width=1\linewidth]{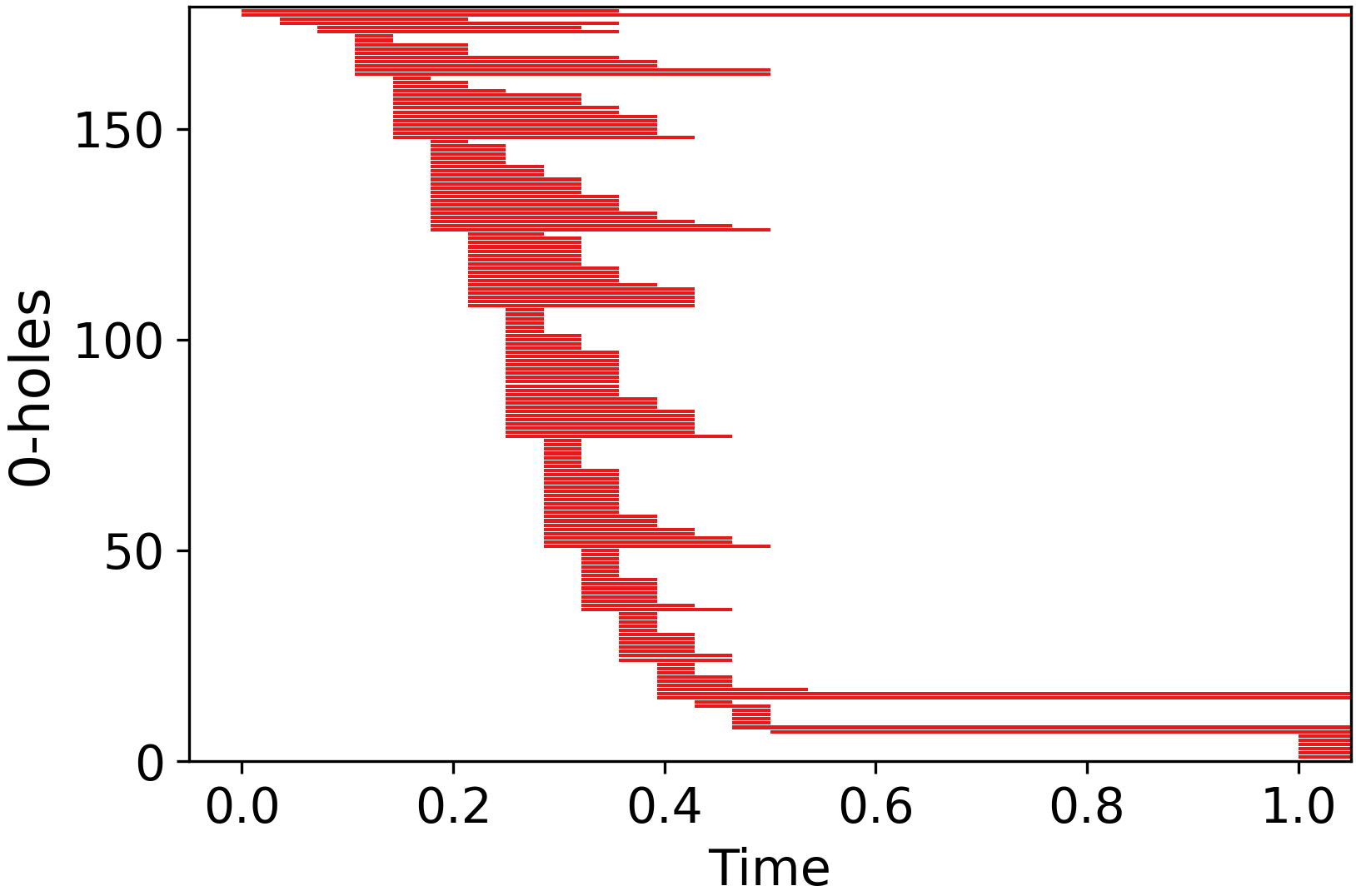} 
		\caption{$\mathcal{F}_\triangle$-induced filtration.} 
		\label{fig:barcodesCube_b}
		\vspace{4ex}
	\end{subfigure} 
	\begin{subfigure}[b]{0.5\linewidth}
		\centering
		\includegraphics[width=1\linewidth]{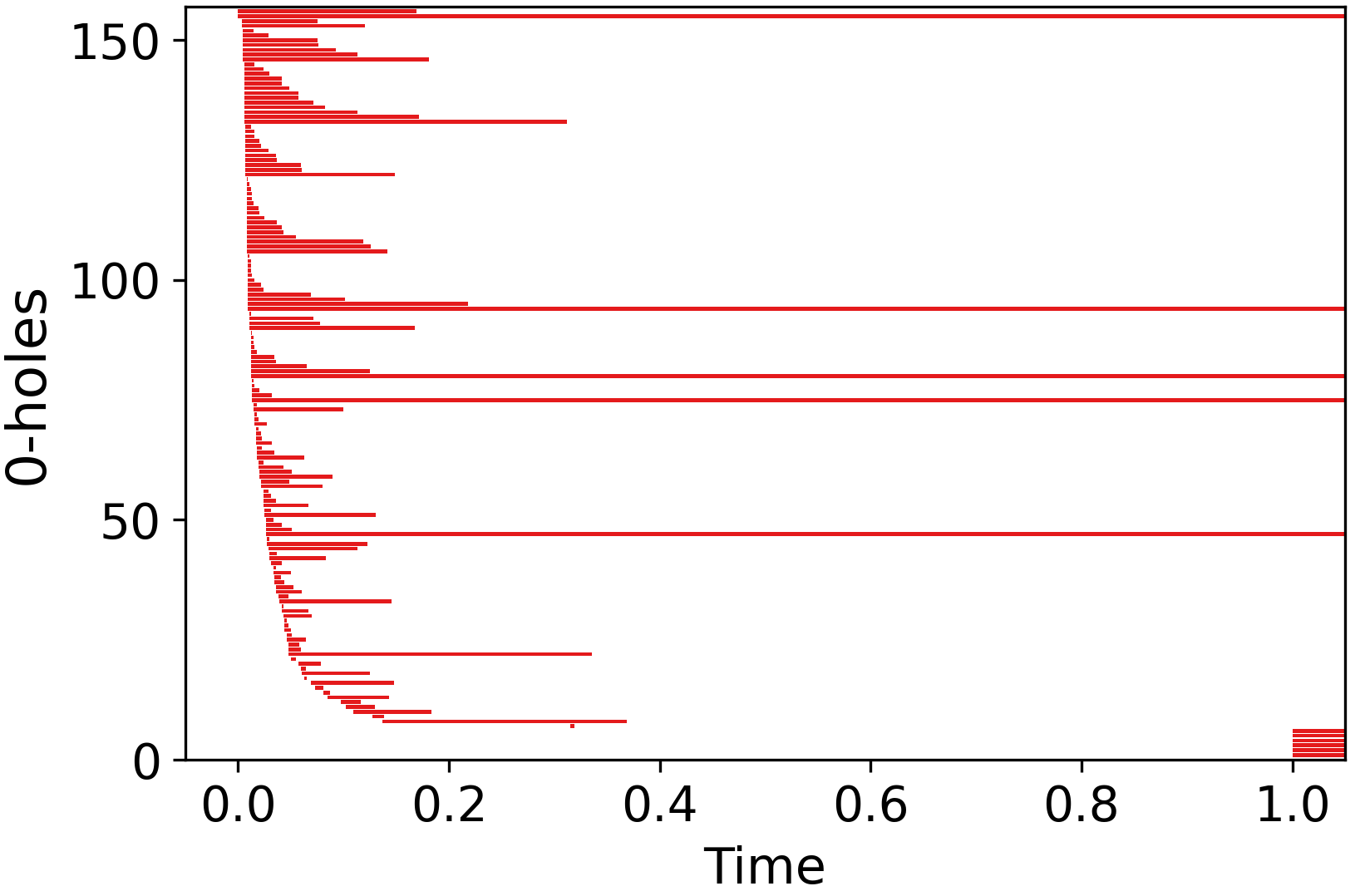} 
		\caption{$\mathcal{F}_{\pentagon}$-induced filtration.} 
		\label{fig:barcodesCube_c}
	\end{subfigure}%% 
	\caption{0-barcodes of a Cube model simulation with average degree 6.}
	\label{fig:barcodesCube}
\end{figure}

\begin{figure*}[h!]
	\centering
	\begin{subfigure}[b]{\linewidth}
		\centering
		\includegraphics[width=\linewidth]{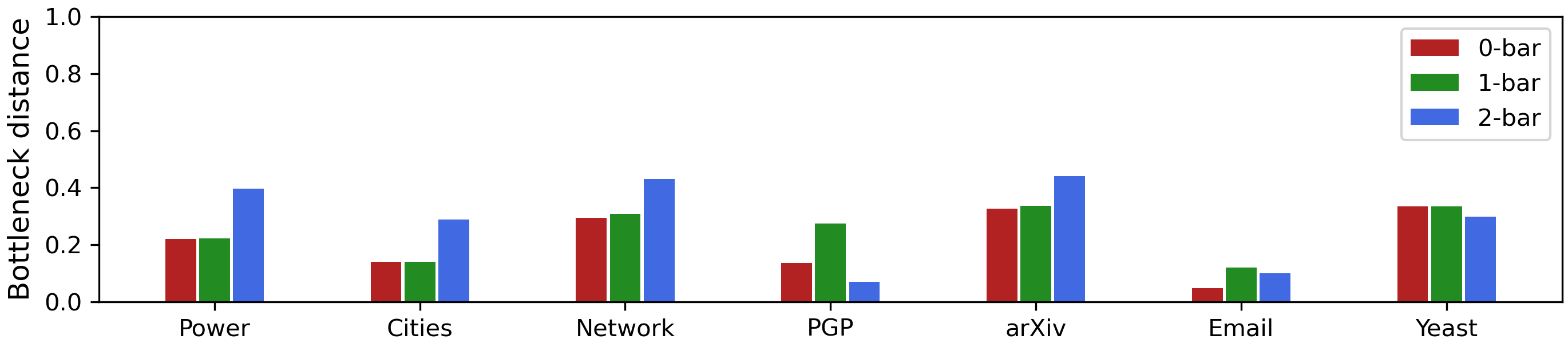} 
		\caption{\small Between $\mathcal{F}$- and $\mathcal{F}_{\triangle}$-induced filtrations.}
		\label{fig:r_plain}
	\end{subfigure}\vspace{2em}
	\begin{subfigure}[b]{\linewidth}
		\centering
		\includegraphics[width=\linewidth]{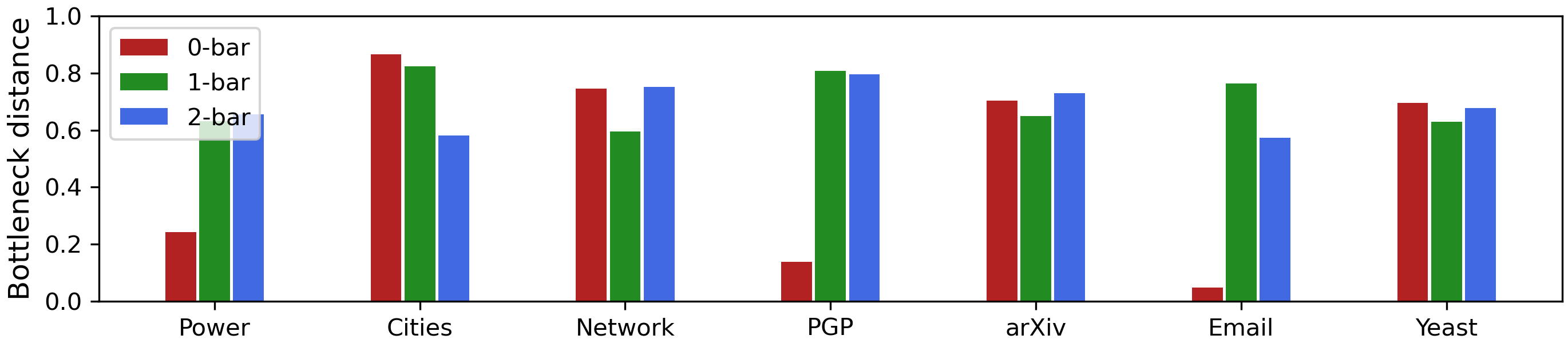}
		\caption{\small Between $\mathcal{F}_\triangle$- and $\mathcal{F}_{\pentagon}$-induced filtrations.} 
		\label{fig:r_pent}
	\end{subfigure}%% 
	\caption{\small Bottleneck distances for the real-world networks.}
	\label{fig:r}
\end{figure*}

We present now the results of our studies; all barcode distances are plotted in Figures \ref{fig:m_plain} to \ref{fig:r}. Let us begin with a general discussion, and then focus on the specifics of each network. We have identified two main reasons for the difference between barcodes:

The first cause, which we will call the ``scale'' effect, is present in most networks and related to the value range of each curvature. Edges whose vertices have low degree are expected to be contained in few cycles, and so their three curvature values will be similar. However, edges between high-degree nodes can belong to a great number of cycles, so their curvature can significantly increase when switching from $\mathcal{F}$ to $\mathcal{F}_\triangle$ and to $\mathcal{F}_{\pentagon}$. Thus, the upper bound of the curvature range increases as we add more cycles, and so normalising the values to $[0,1]$ makes the small ones shift closer to 0. This effect can be observed in Figures \ref{fig:barcodesCube_b} and \ref{fig:barcodesCube_c}, where the bars clearly shorten and shift to the left.

The second cause, which we have named ``the switch'', is related to the signs of the terms of [\ref{eq:4}]. As we have just argued, with enough cycle density the edges with very negative plain curvature (that is, with high-degree nodes) will have highly positive augmented curvature, whereas edges with plain curvature close to 4 will not experiment such drastic changes. Thus, appearance times will be opposite in plain and augmented curvature: edges appearing early (low curvature) with the first will appear late (high curvature) with the second, and vice versa.

We believe the switch has different consequences for barcodes of different dimension. For the 0-barcode, where bars are connected components, we expect many more bars to appear after the switch. With plain curvature the first appearing edges are those with the highest-degree nodes, so the network, which appears from the densest to the sparsest areas, is likely to be pretty well connected at all times. Conversely, with augmented curvature we build the network from sparsest to densest, and so we can expect many more components to appear initially and fuse over time. This effect can be appreciated in Figures \ref{fig:barcodesCube_a} and \ref{fig:barcodesCube_b}.

The consequences of the switch are harder to grasp for higher dimensional barcodes. For simplicity we will focus on 1-holes, but the issues extend directly to higher dimensions. A 1-hole is essentially a cycle, so its appearance time is the highest among the times of its edges. To predict its birth time after the switch we thus need to understand the changes experimented by all of its edges, as their appearance order is likely to have shifted after the switch. Not only that, but also there is not a bijection between cycles and 1-holes: multiple cycles correspond to the same 1-hole, so the appearance time of the latter is the minimum among all generating cycles. This "minimum of a maximum" situation hinders the understanding of 1-hole appearance times. %However, we have observed a pattern that could be interesting to examine more carefully in future work: in some instances the barcodes before and after the switch seem to be, barring some scaling, a negative image of eachother  (see Figure REFERENCIA for a couple of examples, NETWORK 2 Y BA 1). It looks like the appearance time of all the holes has inverted the same way that the appearance time of individual edges would. We believe this explanation to be a bit too naive, and we have not been able to line it up with the theoretically complex appearance time of the bars.

Having discussed the general observed behaviour of the filtrations, let us focus now on each of the studied networks.

\subsection*{Model networks}

\begin{figure}[h!]
	\vspace{3.5em}
	\centering
	\begin{subfigure}[b]{\linewidth}
		\centering
		\includegraphics[width=1\linewidth]{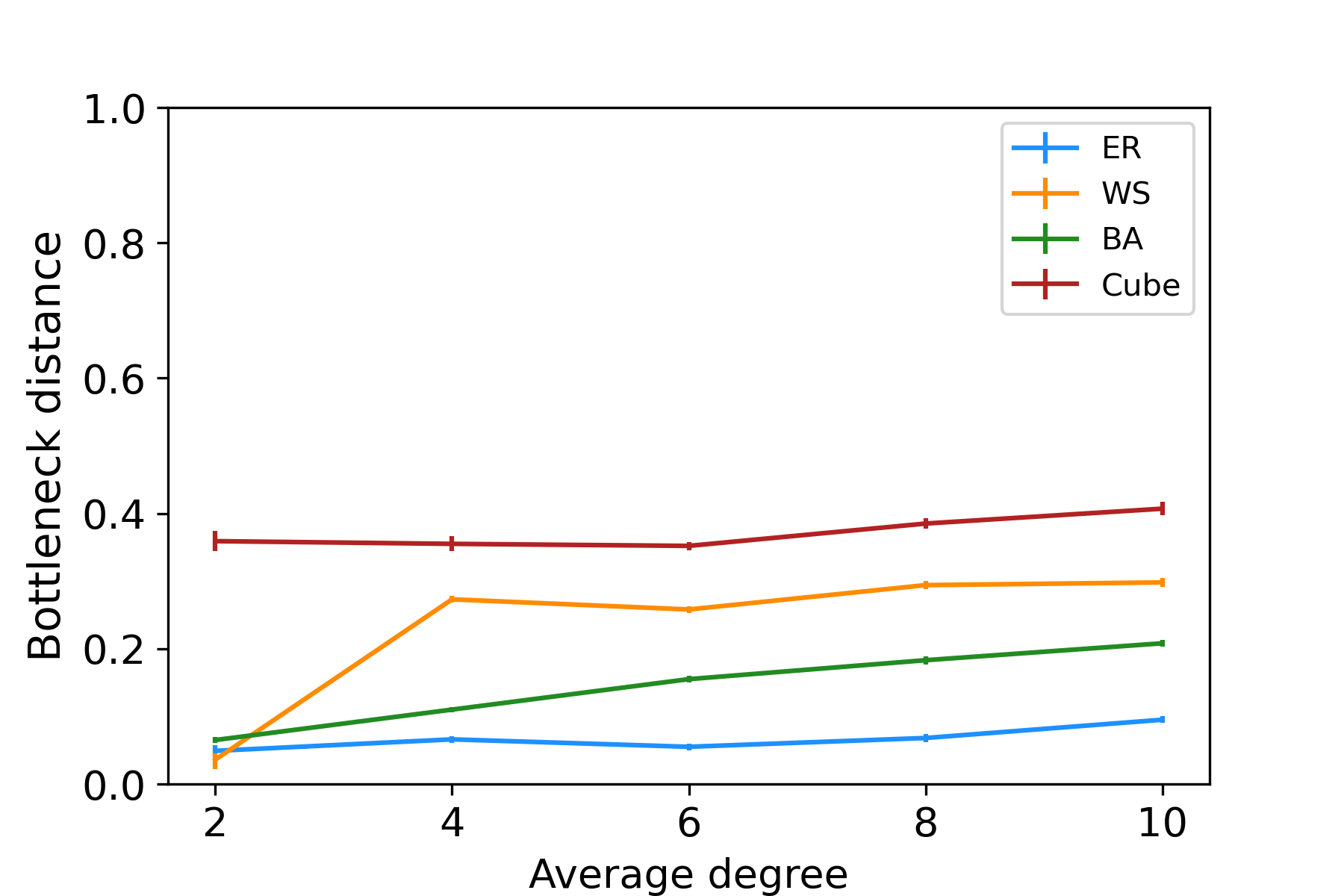}
		\caption{\small 0-barcodes.}
	\end{subfigure}
	\begin{subfigure}[b]{\linewidth}
		\centering
		\includegraphics[width=1\linewidth]{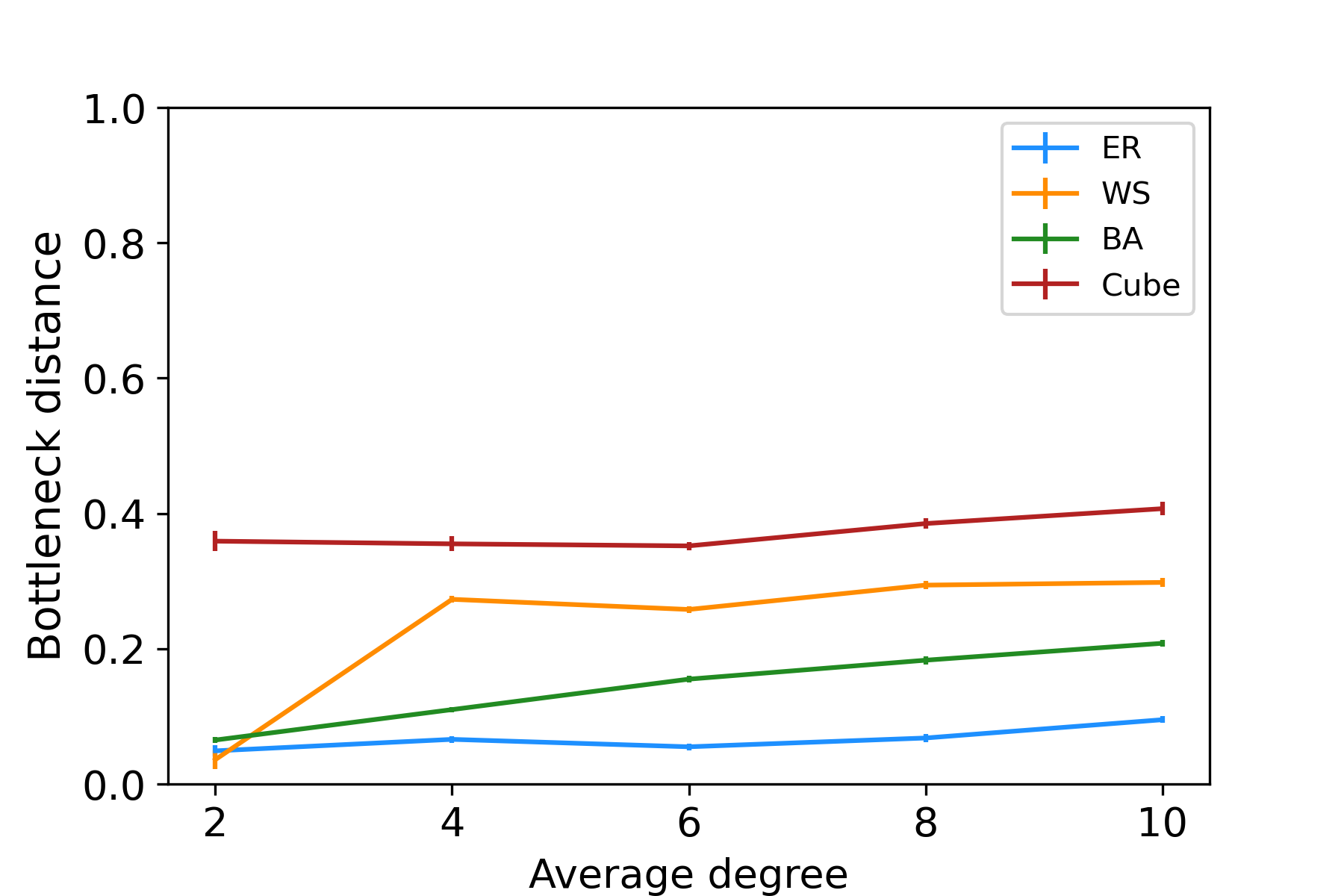} 
		\caption{\small 1-barcodes.} 
	\end{subfigure} 
	\begin{subfigure}[b]{\linewidth}
		\centering
		\includegraphics[width=1\linewidth]{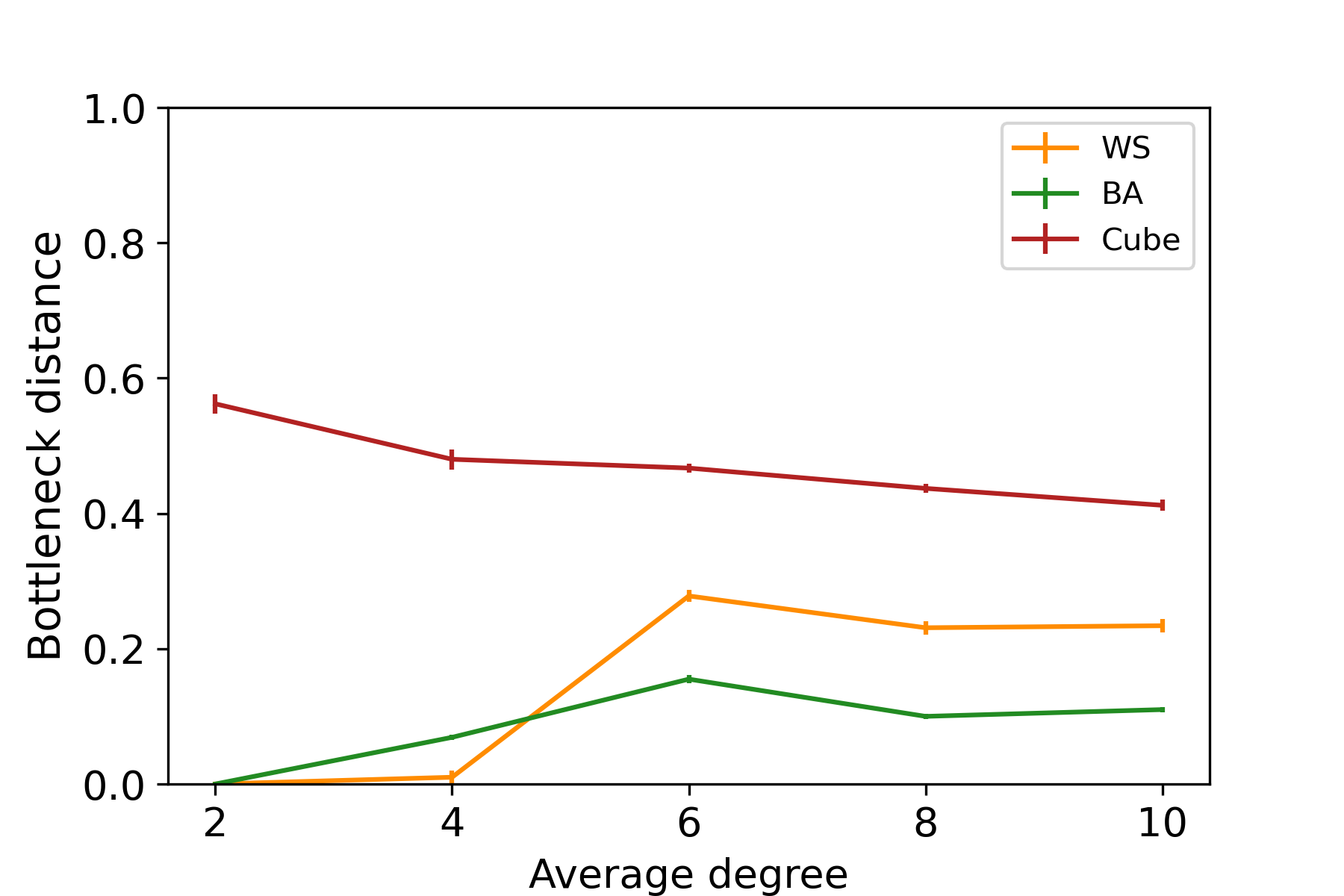} 
		\caption{\small 2-barcodes.} 
	\end{subfigure}%% 
	\caption{\small Average bottleneck distances between the $\mathcal{F}$- and $\mathcal{F}_\triangle$-induced filtrations. Error lines represent the standard error of the mean.}
	\label{fig:m_plain}
	\vspace{3.5em}
\end{figure}

\begin{figure}[h!]
	\vspace{3.5em}
	\centering
	\begin{subfigure}[b]{\linewidth}
		\centering
		\includegraphics[width=1\linewidth]{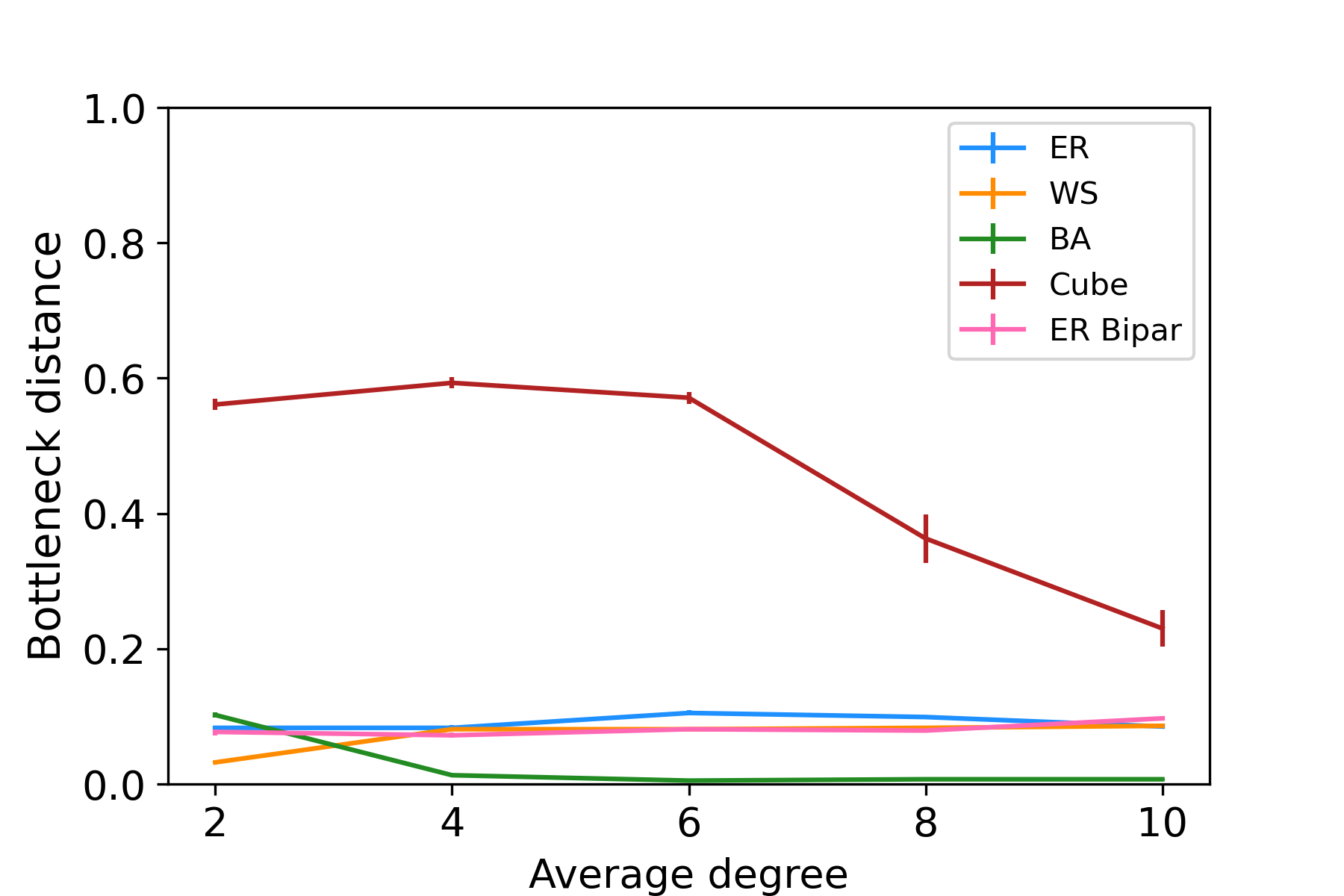}
		\caption{\small 0-barcodes.}
	\end{subfigure}
	\begin{subfigure}[b]{\linewidth}
		\centering
		\includegraphics[width=1\linewidth]{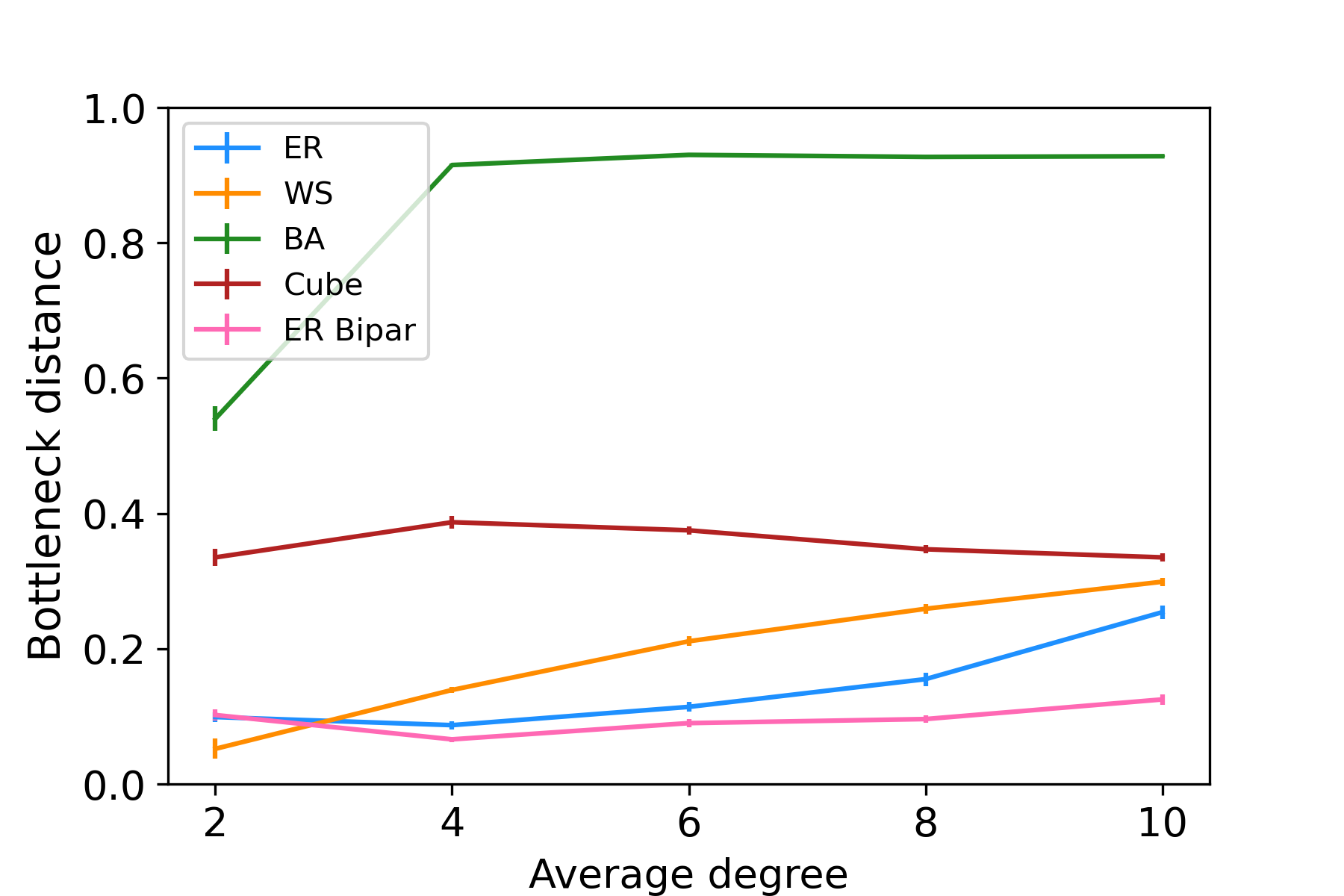} 
		\caption{\small 1-barcodes.} 
	\end{subfigure} 
	\begin{subfigure}[b]{\linewidth}
		\centering
		\includegraphics[width=1\linewidth]{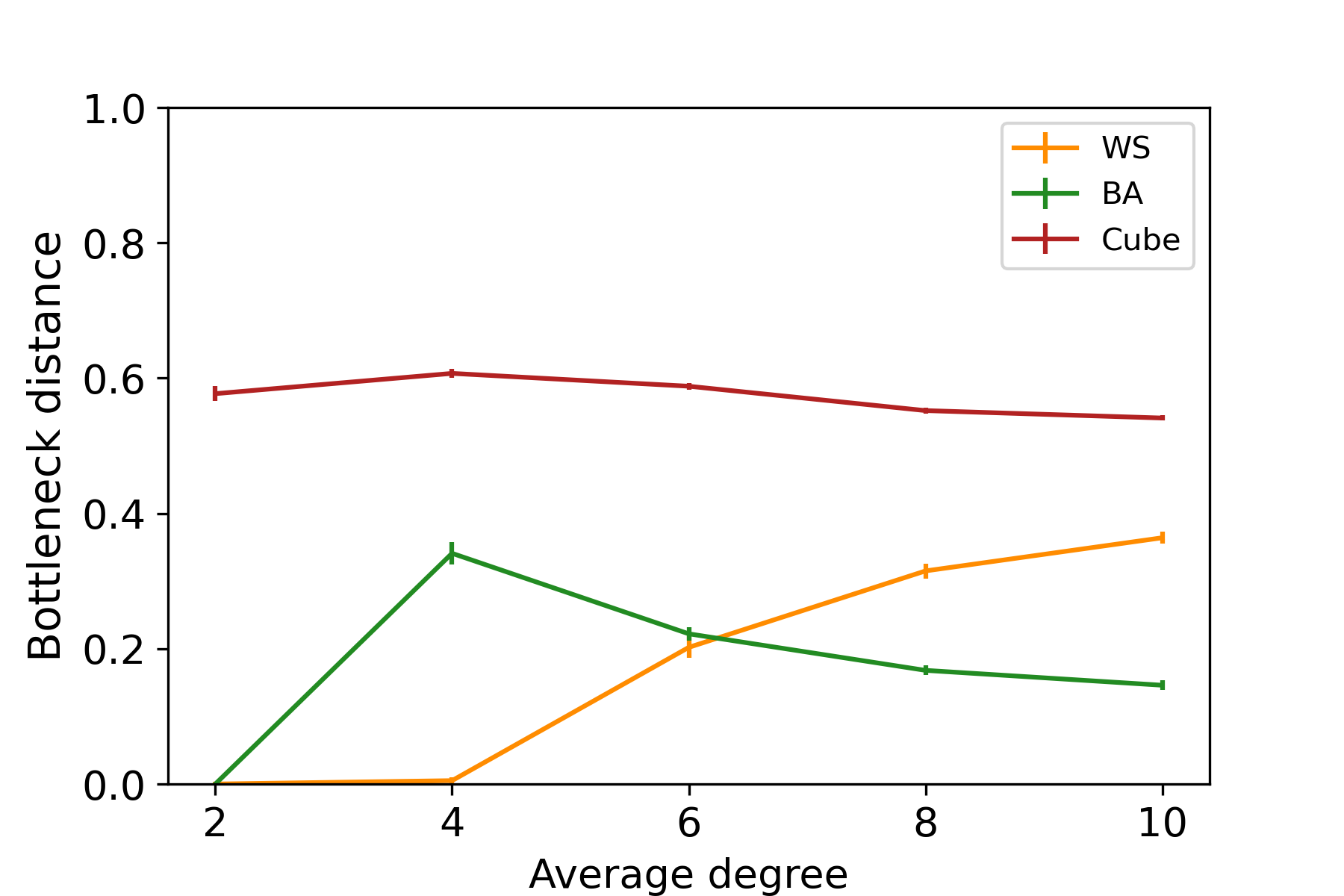} 
		\caption{\small 2-barcodes.} 
	\end{subfigure}%% 
	\caption{\small Average bottleneck distances between the $\mathcal{F}_\triangle$- and $\mathcal{F}_{\pentagon}$-induced filtrations. Error lines represent the standard error of the mean.}
	\label{fig:m_pent}
	\vspace{3.5em}
\end{figure}

Figure \ref{fig:m_plain} shows the average bottleneck distance between the $\mathcal{F}$- and $\mathcal{F}_\triangle$-induced filtrations. Similarly, Figure \ref{fig:m_pent} shows the distances between the $\mathcal{F}_\triangle$- and $\mathcal{F}_{\pentagon}$-induced filtrations. Note that the ER bipartite graph is not present in Figure \ref{fig:m_plain}, as it has no triangles and so $\mathcal{F}$ and $\mathcal{F}_\triangle$ are the same. Note too that no ER data appears for the 2-barcodes. This is because the graphs have a triangle density too low for the triangular faces to clump and form any 2-holes, so there is no 2nd homology group to study. Similarly, the  BA and WS models of average degree 2 have a 2-bottleneck distance of 0 as not enough triangles exist to create a 2-hole.

Looking at Figure \ref{fig:m_plain}, we find that three models exhibit a similar behaviour (ER, WS and BA), whereas the fourth (Cube) shows some noticeable differences. The first three have a low expected number of triangles, which just causes a low-impact, partial switch. The general increase of the distance with the average degree is most probably due to the number of triangles increasing with the degree. One more note about these first three models: the 0-distance for the BA model is almost 0. Plotting the barcodes we observed that they were almost always formed by a single bar born at time 0 -- i.e. at any time step there is only one connected component --, which was already noticed in \cite{FORMANRICCIANDPERSISTENCE}. This is just an extreme case of how with plain curvature (or augmented with few cycles) the network is built from densest to sparsest, as by construction the BA model is very likely to have a small set of nodes connected to all other nodes of the network.

On the other hand, we believe that the big bottleneck distance for the Cube model is a result of its metric building process: let $e=(v_1,v_2)$ be an edge -- that is, $d(v_1,v_2) \leq r$. If we take any other node $v$ connected to $v_1$, there is a probability bigger than $5/16$ (which is the volume proportion of the intersection of two spheres whose centres lie in each other's boundary) that it is also connected to $v_2$. Thus, edges with high-degree nodes are contained in a proportionally high number of triangles, which causes a huge curvature switch (note that $5/16$ is not a tight lower bound, and each triangle increases $\mathcal{F}_\triangle$ by 3). We believe that the decrease of the 0-distance after average degree 6 is due to the final number of connected components. These seem responsible for the $\sim\hspace{-0.2em}0.4$ bottleneck distance, as we usually find multiple infinite bars born at the final steps of the $\mathcal{F}$-filtration that, due to the switch, appear at halfway through the $\mathcal{F}_\triangle$-filtration (see Figures \ref{fig:barcodesCube_a} and \ref{fig:barcodesCube_b}). At average degree 8 the model starts to have one single connected component, and so the bottleneck distance, which then only depends on the intermediate, shorter bars, decreases.

Overall, it seems that results are similar in the trivial case where the models have few triangles, which may be an indication that for some networks we need bigger faces to get a good approximation of the actual value of FR curvature. On the other hand, we have found important differences on the one model with a substantial amount of triangles (Cube), so it seems that we should not use plain FR curvature for these cases, as one would expect.

As per Figure \ref{fig:m_pent}, it seems like the addition of squares and pentagons is enough to cause a significant curvature switch on the ER (standard and bipartite), WS, and BA models. The bottleneck distances of the first three still are not too high, but the $\mathcal{F}_{\pentagon}$-induced 0-barcodes generally show many more, much shorter, left-shifted bars (similar to Figure \ref{fig:barcodesCube_c}), which we believe are fruit of the curvature switch and scaling effect.

\begin{figure}%[tbhp]
	\centering
	\begin{subfigure}[b]{0.47\linewidth}
		\centering
		\includegraphics[width=1\linewidth]{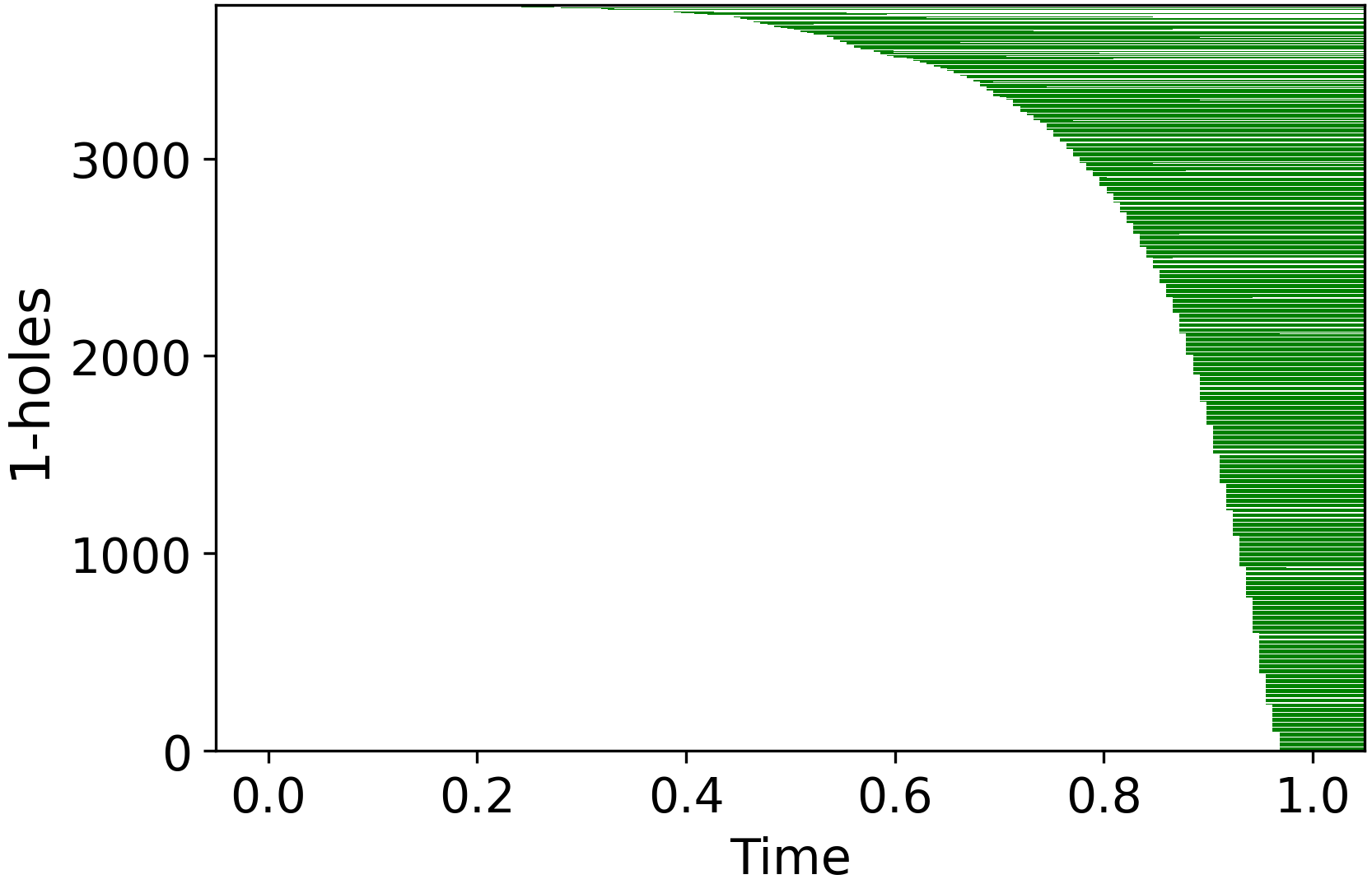} 
		\caption{1-barcode, $\mathcal{F}_\triangle$-induced.} 
		\label{fig:barcodesBA_a}
		\vspace{4ex}
	\end{subfigure}%
	\hspace{1em}%
	\begin{subfigure}[b]{0.47\linewidth}
		\centering
		\includegraphics[width=1\linewidth]{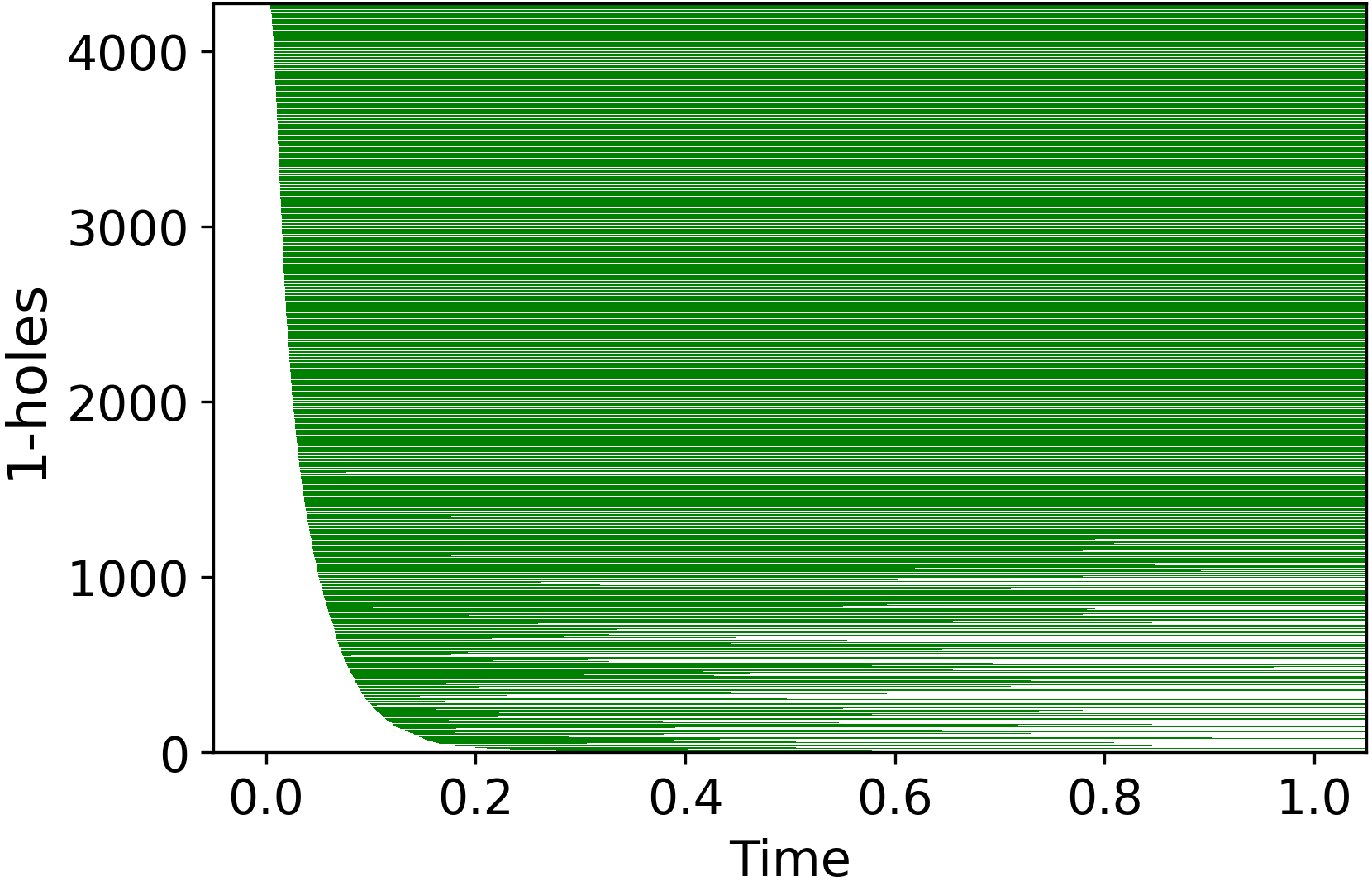} 
		\caption{1-barcode, $\mathcal{F}_{\pentagon}$-induced.} 
		\label{fig:barcodesBA_b}
		\vspace{4ex}
	\end{subfigure} 
	\begin{subfigure}[b]{0.47\linewidth}
		\centering
		\includegraphics[width=1\linewidth]{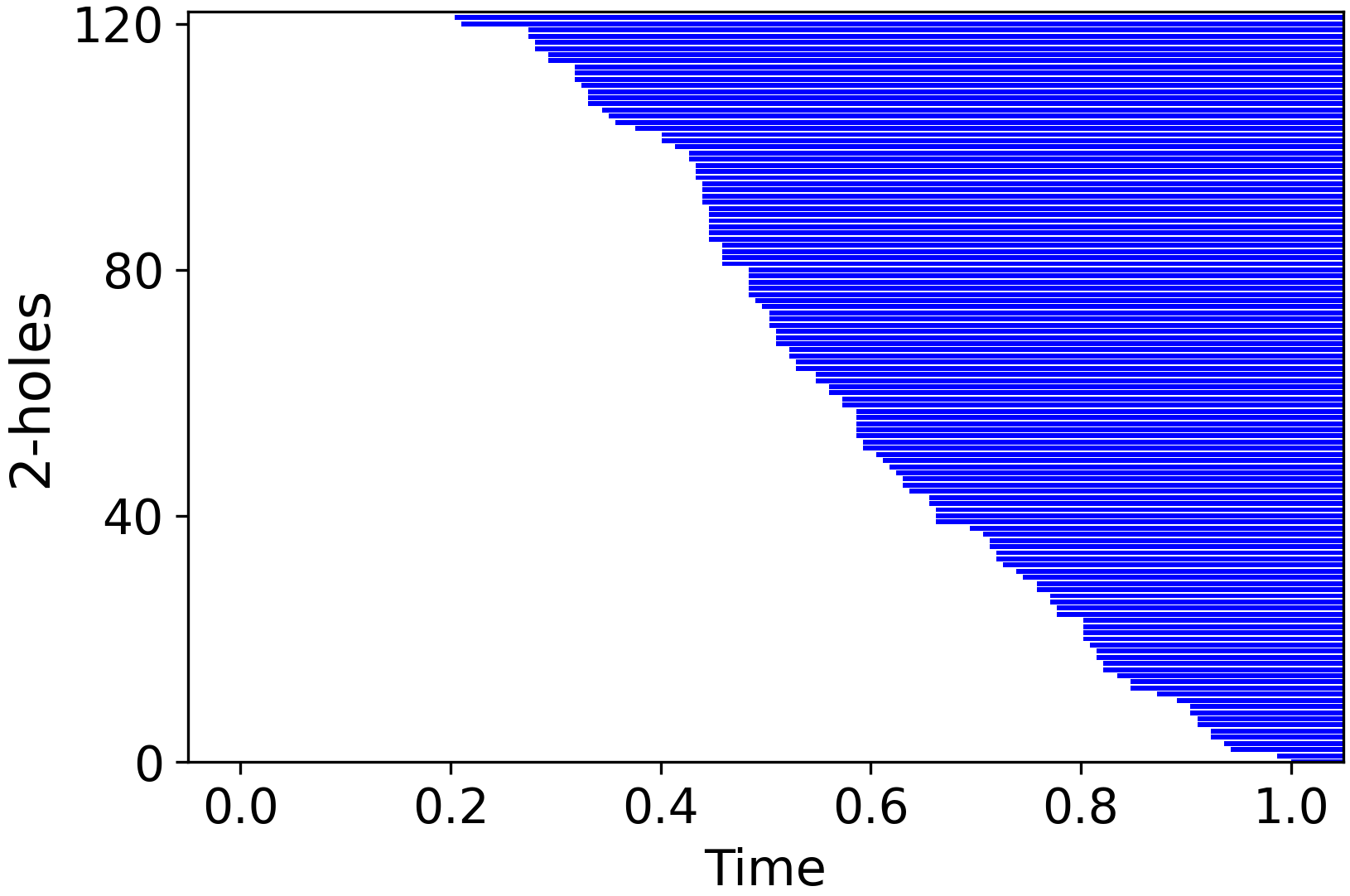} 
		\caption{2-barcode, $\mathcal{F}_\triangle$-induced.} 
		\label{fig:barcodesBA_c}
	\end{subfigure}%
	\hspace{1em}%
	\begin{subfigure}[b]{0.47\linewidth}
		\centering
		\includegraphics[width=1\linewidth]{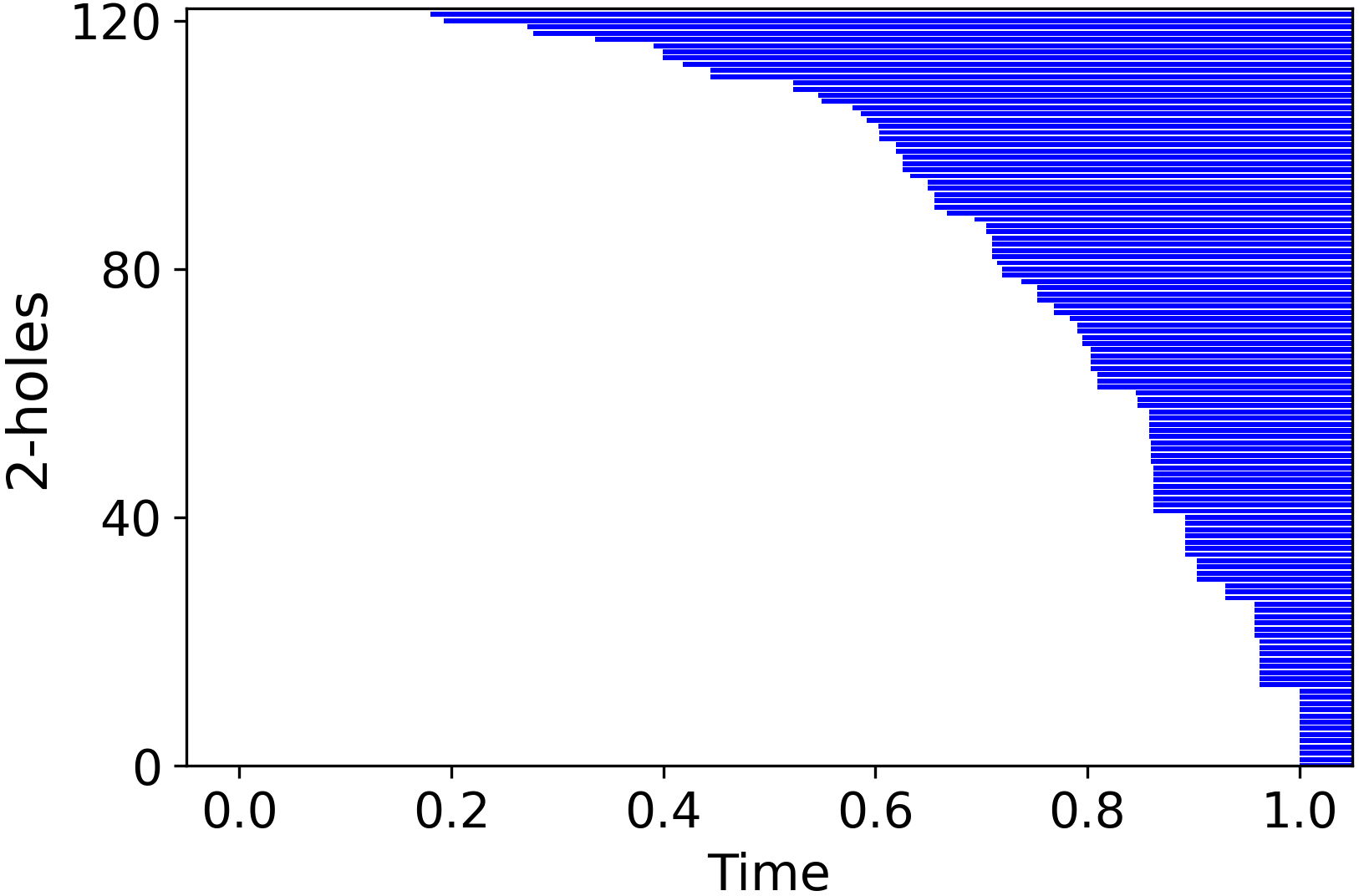} 
		\caption{2-barcode, $\mathcal{F}_{\pentagon}$-induced.} 
		\label{fig:barcodesBA_d}
	\end{subfigure} %% 
	\caption{Some barcodes of a BA model simulation with average degree 6.}
	\label{fig:barcodesBA}
\end{figure}

Regarding the BA model, there is a great distance between its 1-barcodes (see some examples in Figure \ref{fig:barcodesBA_a} and \ref{fig:barcodesBA_b}), which we attribute to the building process generating a huge amount of squares and a massive number of pentagons \cite{NUMBEROFLOOPSBA}, which causes a full switch. We do not know whether the 1-barcode difference is the result of just a very strong scaling effect, or there is an additional phenomenon happening -- as the 2-barcodes do not seem to be so affected by curvature scaling, see Figures \ref{fig:barcodesBA_c} and \ref{fig:barcodesBA_d}. An explanation for the invariance of the 2-barcode could be the fact that there are many edges involved in the formation of a 2-hole (the smallest possible is an empty tetrahedron), and the maximum appearance time of a set of edges can potentially remain invariant under a curvature switch if the set is big enough. Finally, the 0-barcode distance is again almost 0, but this time as a result of curvature scaling: the switch causes a lot of short bars to initially appear in the $\mathcal{F}_{\pentagon}$-filtration, but the curvature scaling shortens them so much that they barely affect the barcode distance.

As per the Cube model, the curvature switch has already occurred in the $\mathcal{F}_\triangle$-filtration, and the barcode distances are mostly due to curvature scaling, which shifts all bars to the left. The drop on the 0-distance at average degree 6 follows the same explanation as before -- now the initial distance of $\sim\hspace{-0.2em}0.6$ is a result of the infinite bars shifting from being born halfway through the filtration to being born during its first moments (see Figures \ref{fig:barcodesCube_b} and \ref{fig:barcodesCube_c}). 

Overall, although adding squares and pentagons seems to cause a curvature switch in the models which had few triangles, we fear that not taking orientations into account (which would make some cycles cancel their contributions out) may create a scaling effect that heavily distorts the persistence diagrams. We have observed how the appearance times of most bars end up very close to $0$, and non-infinite bars are drastically shortened. We believe that the combination of these two effects makes the bottleneck distance quite inconsequential, and causes barcodes as a whole to lose a substantial amount of information.

\subsection*{Real-world networks} Figures \ref{fig:r_plain} and \ref{fig:r_pent} show the results of our studies. There does not seem to be a clear relationship between the average degrees and the bottleneck distances, which should not be surprising -- we have already seen (e.g. Cube model) that other structural characteristics also play a big part in the results.

It seems like the overall impact of the squares and pentagons is much bigger than that of the triangles. However, examination of all the individual barcodes suggests that this is a result of curvature scaling (similarly to what happens to low degree Cube models) rather than the triangles being insufficient to trigger the curvature switch.

Regarding the difference between the $\mathcal{F}$- and $\mathcal{F}_\triangle$-induced filtrations, Figure \ref{fig:r_plain}, most 0-barcodes exhibit a behaviour that would be expected from a curvature switch: the number of 0-bars grows between a $25\%$ and a $200\%$ for most networks, and the infinite 0-bars experiment a shift to the left. Still, a clearer indication that triangles cannot generally be omitted can be found in the 1-barcodes, where we have observed an effect that did not appear in any of the models: some of the 1-barcodes see a sharp increase in the number of bars. In particular, for the Network network the number of 1-bars goes from 12 to 191; for PGP, from 1315 to 2679; for arXiv, from 910 to 2069; and for Yeast, from 84 to 1725. Even if we cannot fully explain why this is the case (as we lack insight of both the networks' structure and the behaviour of 1-holes), we believe it is strong evidence against using plain FR curvature as a substitute of augmented curvature, as their results greatly differ.

As per the difference between the $\mathcal{F}_\triangle$- and $\mathcal{F}_{\pentagon}$-induced filtrations, Figure \ref{fig:r_pent}, we believe that all distances are a result of an extreme curvature scaling effect -- at least for the 0- and 1-barcodes. The 0- and 1-bars of all networks are drastically shifted to the left, and most of the non-infinite bars are severely shortened -- a behaviour similar to Figures \ref{fig:barcodesCube_c} and \ref{fig:barcodesBA_b}. We have observed that most barcodes have infinite bars born in the middle of the filtration, which, just like with the Cube model, shift to an appearance time of $\sim \hspace{-0.2em}0$ with $\mathcal{F}_{\pentagon}$, which causes the high bottleneck distances. Regarding the few barcodes for which this is not the case (0-barcodes of Power, PGP, and Email), their behaviour is similar to the BA model: all of them have a single infinite bar, and thus the bottleneck distance is just due to the mismatch of the shorter, intermediate bars that are severely reduced with $\mathcal{F}_{\pentagon}$. Interestingly, similarly to the BA model (Figure \ref{fig:barcodesBA_d}), most of the $\mathcal{F}_{\pentagon}$-induced 2-barcodes do not show such an extreme scaling effect.

As with the model networks, the results seem to suggest that using $\mathcal{F}_{\pentagon}$ for non-quasiconvex networks is too rough of an approximation, as the fact that [\ref{eq:4}] does not take into account the orientation of the faces causes edges with high-degree nodes to have disproportionately big curvatures.

\section*{Conclusion}
We have provided a short overview of persistent homology and FR curvature, and discussed how to compute augmented FR curvature for non-quasiconvex augmented networks. The focus of our work has then been analysing the technique proposed in \cite{FORMANRICCIANDPERSISTENCE} to use FR curvature to build time filtrations for the persistent homology study of networks. We have used three different versions of FR curvature: plain, triangle-augmented (in a version theoretically adequate for all networks), and pentagon-augmented (in a version theoretically adequate only for quasiconvex networks). We have studied both model networks and real-world networks with different structures and average degrees.

Our results suggest that using plain curvature as a computationally faster version of triangle-augmented curvature may not produce accurate results. The only networks in which the omission of the triangular faces did not suppose a big difference of results was the trivial case in which networks did not have many triangles to begin with.

On the other hand, our study seems to also advise against using $\mathcal{F}_{\pentagon}$ for non-quasiconvex augmented networks. For these, the formula seems to give excessive weight to the 4- and 5-cycles, as it does not take into account their orientations. Not only does this cause an offset from the actual curvature value, but also leads to a severe curvature scaling effect that renders bottleneck distance meaningless -- as most bars are drastically shortened and shifted to a low appearance time.

Still, the results indicate that for networks with a low number of triangles (e.g BA model) one may need to take into account higher cycles in order to have an accurate approximation of augmented FR curvature. Nonetheless, we believe that one should employ the triangle augmentation until the general, augmented FR curvature (formula [\ref{eq:6}]) is further researched and implemented, as we think that the application of $\mathcal{F}_{\pentagon}$ to general networks leads to heavily distorted results.

A possible avenue of future research would be further studying [\ref{eq:6}], both theoretically and from a computational implementation point of view. We believe this would be fruitful for all Forman-Ricci curvature studies, not only persistent homology ones. Another approach could be studying the behaviour of the appearance times of $k$-holes for $k\geq 1$, which would provide insight on the nature and interpretation of such holes in networks.

\acknow{This paper started as a final project for the Networks module of Oxford's MSc in Mathematical Sciences. We would like to thank Prof. Lambiotte, as well as class tutors Ms. Semenova and Mr. Falco, for their guidance throughout the course.
	
This article has been created using PNAS Tex article template.}

\showacknow{}

% Bibliography
\bibliography{pnas-sample}

% Appendix
%\onecolumn
%\section*{Appendix}

\end{document}